\documentclass[a4paper]{article}
\usepackage[margin=1in]{geometry}
\usepackage{cite}
\usepackage[T1]{fontenc}
\usepackage[utf8]{inputenc}
\usepackage{amsmath}
\usepackage{amsfonts}
\usepackage{amssymb}
\usepackage{amsthm}
\usepackage{graphicx}

\usepackage[capitalize]{cleveref}
\newcommand{\adj}{\textup{H}}
\newcommand{\mA}{\textup{A}}
\newcommand{\mB}{\textup{B}}
\newcommand{\xC}{\mathbb{C}}
\newcommand{\mG}{\textup{G}}
\newcommand{\vn}{\mathbf{n}}
\newcommand{\vq}{\mathbf{q}}
\newcommand{\xR}{\mathbb{R}}
\newcommand{\vr}{\mathbf{r}}
\newcommand{\vs}{\mathbf{s}}
\newcommand{\mT}{\textup{T}}
\newcommand{\vu}{\mathbf{u}}
\newcommand{\mU}{\textup{U}}
\newcommand{\vv}{\mathbf{v}}
\newcommand{\mV}{\textup{V}}
\newcommand{\vw}{\mathbf{w}}

\newcommand{\vzero}{\mathbf{0}}
\newcommand{\vone}{\mathbf{1}}

\newcommand{\aabs}[1]{|#1|}
\newcommand{\norm}[1]{\left\|#1\right\|}
\newcommand{\nnorm}[1]{\|#1\|}
\DeclareMathOperator*{\argmax}{arg\,max}
\DeclareMathOperator*{\argmin}{arg\,min}

\newtheorem{remark}{Remark}
\crefname{remark}{Remark}{Remarks}
\newtheorem{theorem}{Theorem}
\crefname{theorem}{Theorem}{Theorems}
\crefformat{equation}{\textup{(#2#1#3)}}

\usepackage{tikz}
\usepackage{pgfplots}
\usepackage{pgfplotstable}
\pgfplotsset{compat=1.9}

\begin{document}
\title{Adaptive approximation of nonlinear eigenproblems by minimal rational interpolation}

\author{Davide \textsc{Pradovera}\thanks{University of Vienna, Oskar-Morgenstern-Platz 1, 1090 Vienna, Austria. E-mail: \texttt{davide.pradovera@univie.ac.at}}} 

\maketitle

\begin{abstract}
\noindent
We describe a strategy for solving nonlinear eigenproblems numerically. Our approach is based on the approximation of a vector-valued function, defined as solution of a non-homogeneous version of the eigenproblem. This approximation step is carried out via the \emph{minimal rational interpolation} method. Notably, an adaptive sampling approach is employed: the expensive data needed for the approximation is gathered at locations that are optimally chosen by following a greedy error indicator. This allows the algorithm to employ computational resources only where ``most of the information'' on not-yet-approximated eigenvalues can be found. Then, through a post-processing of the surrogate, the sought-after eigenvalues and eigenvectors are recovered. Numerical examples are used to showcase the effectiveness of the method.
\end{abstract}

\section{Introduction}\label{sec:1}
Let $\mT:\xC\to\xC^{n\times n}$ be a matrix-valued function, with $n$ a (large) integer. We are interested in identifying \emph{eigenpairs} $(\lambda,\vw_\lambda)\in\xC\times\xC^n$ such that
\begin{equation}\label{eq:ep}
\mT(\lambda)\vw_\lambda=\vzero,\quad\text{with }\vw_\lambda\neq\vzero.
\end{equation}
More precisely, we seek \emph{all} the eigenpairs whose eigenvalues $\lambda$ lie in a certain region $Z\subset\xC$ of the complex plane, e.g., in a half-plane, on a line segment, or inside a Bernstein ellipse. Problems of this kind appear in several branches of science and engineering. For instance, $\lambda$ might denote the frequency of vibration when studying the time-harmonic behavior of mechanical or electrical systems. In such cases, the eigenvector $\vw_\lambda$ provides information on the resonating modes of the system at hand.

Many algorithms have been proposed for finding eigenpairs satisfying \cref{eq:ep}. Among these, we mention methods based on (affine or rational) approximations of $T(\cdot)$, which then cast an approximated version of \cref{eq:ep} as an augmented \emph{linear} generalized eigenproblem $\mA\widehat{\vw}_{\widehat{\lambda}}=\widehat{\lambda}\mB\widehat{\vw}_{\widehat{\lambda}}$. See, e.g., \cite{guttel}. One of the main drawbacks of this class of methods is that the linearized eigenproblem might be rather large, especially when the matrix-valued function $\mT$ is ``hard to linearize'', i.e., when many terms are necessary in its affine/rational approximation. On the other hand, we can find approaches based on integration of $\mT(\cdot)^{-1}$, or of $\mT(\cdot)^{-1}\mV$, with $\mV\in\xC^{n\times m}$ a tall and skinny matrix ($m\ll n$), over the contour $\partial Z$. See, e.g., \cite{beyn}. Such methods rely on the fact that, under standard assumptions (see, e.g., \cref{th:keldysh}), eigenvalues of \cref{eq:ep} coincide with poles of $\mT(\cdot)^{-1}$. Here, one of the main difficulties is finding a reliable quadrature formula to approximate the necessary contour integrals. Notably, one should strive to reduce the number of quadrature points, since an expensive evaluation of $\mT(\cdot)^{-1}\mV$ (i.e., the assembly of $\mT(\cdot)$ and the solution of $m$ linear systems of size $n$) is needed at each of them.

In this work, we follow a logic similar to that used in contour-based methods: given a vector $\vv\in\xC^n$, we consider the vector-valued function\footnote{We can replace the single vector $\vv$ with a tall and skinny matrix $\mV\in\xC^{n\times m}$ as discussed for contour-based methods. One could elect to do this, e.g., as a way to reduce the risk of ``missing'' eigenvectors that are (almost) orthogonal to $\vv$, cf.~\cref{th:keldysh}. This makes $\mU:\xC\ni z\mapsto\mT(z)^{-1}\mV\in\xC^{n\times m}$ matrix-valued (but with only few columns). Our discussion still applies, with the only change being the need to replace some Euclidean scalar products with Frobenius ones.} $\vu:\xC\ni z\mapsto\mT(z)^{-1}\vv\in\xC^n$ and approximate it using the \emph{minimal rational interpolation} method \cite{mri}. This yields a vector-valued rational function $\widetilde{\vu}$ that approximates $\vu$ over the region of interest $Z$. As we will see in \cref{th:keldysh}, the poles of $\vu$ can be related to the eigenvalues of \cref{eq:ep}, and the corresponding residues to the eigenvectors. As such, poles and residues of the surrogate $\widetilde{\vu}$, which are much simpler to compute than those of $\vu$ (see \cref{sec:2}), provide approximations of the eigenpairs of interest.

In our approach, the accuracy of the approximation of the eigenpairs is, in a loose sense, inherited from the accuracy with which $\widetilde{\vu}$ approximates $\vu$. For this reason, we aim at obtaining the ``best'' (in some sense to be specified) approximation $\widetilde{\vu}$, given a computational budget. More precisely, we assume that we can only evaluate $\vu$ a fixed number of times.

In order to achieve our claimed ``optimality of approximation'' for a given computational budget, we endeavor to place the sample points (i.e., the locations $z_1,z_2,\ldots\in\xC$ where $\vu$ is to be evaluated) in such a way as to maximize the ``amount of information'' that each sample provides. For this, we use the \emph{adaptive} version of minimal rational interpolation, originally introduced in \cite{zmri}. The method is presented in \cref{sec:2} and, in summary, consists of an incremental strategy for adding new sample points, each chosen in a ``greedy'' way based on the previous ones, until the computational budget is exhausted.

Before proceeding further, we remark that, due to its non-intrusive nature, cf.~\cref{sec:4}, our method may also be applied in settings more general than the one presented above. For instance, we might use it even for nonlinear eigenproblems \emph{with eigenvector nonlinearities}, e.g.,
\begin{equation*}
\mT(\lambda,\vw_\lambda)\vw_\lambda=\vzero,
\end{equation*}
with $\mT:\xC\times\xC^n\to\xC^{n\times n}$, cf. \cite{lietaert}. Of course, this is feasible only under the assumption that an (approximate) solver is available for the \emph{non-homogeneous} nonlinear version of the problem, i.e., 
\begin{equation*}
\text{for given }z\in\xC\text{ and }\vv\in\xC^n,\text{ find }\vu(z)\in\xC^n\text{ such that }\mT(z,\vu(z))\vu(z)=\vv.
\end{equation*}
Otherwise stated, one must be able to evaluate $\vu(z)$ (the solution of the above non-homogeneous problem) at any value of $z$.

\section{Approximation method}\label{sec:2}
We start by describing the strategy for approximating a generic function $\vu:\xC\to\xC^n$, given the samples $\{(z_j,\vu(z_j))\}_{j=1}^S$, with $z_1,\ldots,z_S$ being \emph{given} distinct complex points (although our method can be modified to account for point confluence), using minimal rational interpolation. We only provide the bare algorithm here, referring to \cite{mri,mrib} for more details.

We seek a rational approximation $\widetilde{\vu}:\xC\to\xC^n$ that \emph{interpolates} the data exactly. For this, we rely on the extremely compact \emph{barycentric format}:
\begin{equation}\label{eq:approx}
\widetilde{\vu}(z)=\frac{\sum_{j=1}^Sq_j\vu(z_j)\prod_{j'\neq j}(z-z_{j'})}{\sum_{j=1}^Sq_j\prod_{j'\neq j}(z-z_{j'})}=\underbrace{\sum_{j=1}^S\frac{q_j\vu(z_j)}{z-z_j}}_{\vn(z)}\Bigg/\underbrace{\sum_{j=1}^S\frac{q_j}{z-z_j}}_{d(z)},
\end{equation}
which achieves interpolation of the data for all (nonzero) complex coefficients $q_1,\ldots,q_S$. The \emph{minimal rational interpolant} of the data is defined as a rational function as in \cref{eq:approx} whose barycentric coefficients $q_1,\ldots,q_S$ satisfy
\begin{equation}\label{eq:optimal}
(q_1,\ldots,q_S)=\argmin_{q_1',\ldots,q_S'\in\xC}\frac{\nnorm{\sum_{j=1}^Sq_j'\vu(z_j)}_2^2}{\aabs{\sum_{j=1}^Sq_j'}^2}=\argmin_{\sum_{j=1}^Sq_j'=1}\underbrace{\sum_{j,j'=1}^S\overline{q_{j'}'}q_j'\vu(z_{j'})^\adj\vu(z_j)}_{\nnorm{\sum_{j=1}^Sq_j'\vu(z_j)}_2^2}.
\end{equation}
The last step follows by a simple scaling argument $q_j'\mapsto q_j'\big/\sum_{j'=1}^Sq_{j'}'$. The bar $\overline{\cdot}$ and the subscript $\cdot^\adj$ denote complex conjugation and conjugate transposition, respectively. We remark that, since
\begin{equation*}
\lim_{z\to\infty}\widetilde{\vu}(z)=\frac{\sum_{j=1}^Sq_j\vu(z_j)}{\sum_{j=1}^Sq_j},
\end{equation*}
cf.~\cref{eq:approx}, this can be interpreted as minimizing the value of $\nnorm{\widetilde{\vu}}_2$ at infinity.

The minimization in \cref{eq:optimal} is a quadratic programming problem, whose first-order optimality (KKT) conditions read
\begin{equation*}
\left\{\hspace{-.5em}\begin{array}{l}
\mG\vq=\lambda\vone,\\
\vone^\adj\vq=1,
\end{array}\right.\text{ where }
\vq=\begin{bmatrix}
q_1 \\ \vdots \\ q_S
\end{bmatrix}\text{, }\mG=\begin{bmatrix}
\vu(z_1)^\adj\vu(z_1) & \cdots & \vu(z_1)^\adj\vu(z_S) \\ \vdots & \ddots & \vdots \\ \vu(z_S)^\adj\vu(z_1) & \cdots & \vu(z_S)^\adj\vu(z_S)
\end{bmatrix}\text{, and }\vone=\begin{bmatrix}
1 \\ \vdots \\ 1
\end{bmatrix}\in\xC^S.
\end{equation*}
Note that $\lambda\in\xC$ is a Lagrange multiplier, introduced to enforce the linear constraint. An explicit solution of \cref{eq:optimal} can be obtained by left-multiplying\footnote{Inverting the Gramian matrix $\mG$ may be numerically unstable. A more robust formulation for minimal rational interpolation also exists, where the linear constraint $\sum_jq_j'=1$ in \cref{eq:optimal} is replaced by the Euclidean normalization $\sum_j\big|q_j'\big|^2=1$. In this case, the closed-form solution to the problem is: given the eigendecomposition $\mG=\mV\Lambda\mV^\adj$, with the diagonal elements of $\Lambda$ in non-decreasing order, the optimal $\vq$ is the first column of $\mV$. See \cite{mrib} for more details.} the first KKT condition by $\mG^{-1}$ (thus finding $\vq$ as a function of $\lambda$) and then plugging the corresponding value of $\vq$ in the second KKT condition to find $\lambda$. This yields $\vq=\mG^{-1}\vone~\big/\left(\vone^\adj\mG^{-1}\vone\right)$ as optimal solution.

\subsection{Eigenpair-finding by post-processing of minimal rational interpolant}
Once a minimal rational interpolant has been computed, we obtain the eigenpairs of interest by a pole- and residue-finding step, which we detail here.

Given a rational function $\widetilde{\vu}$ in barycentric form, see \cref{eq:approx}, its poles can be found as the (at most) $S-1$ finite eigenvalues of a $(S+1)\times(S+1)$ generalized eigenproblem in arrowhead form, cf., e.g., \cite{aaa}. More specifically, we have the characterization: for any finite $\lambda\in\xC$,
\begin{equation}\label{eq:poles}
\lim_{z\to\lambda}\norm{\widetilde{\vu}(z)}=\infty\quad\Leftrightarrow\quad\exists\vv\neq\vzero\text{ such that }\begin{bmatrix}
0 & q_1 & \cdots & q_S \\
1 & z_1 & & \\
\vdots & & \ddots & \\
1 & & & z_S \\
\end{bmatrix}\vv=\lambda\begin{bmatrix}
0 & & & \\
 & 1 & & \\
 & & \ddots & \\
 & & & 1 \\
\end{bmatrix}\vv.
\end{equation}

Now, assume that $\lambda\in\xC$ is a pole of $\widetilde{\vu}$ with order $N>0$, i.e., a (finite) eigenvalue of \cref{eq:poles} with multiplicity $N$. We seek the corresponding $N$ residues $\vr_1,\ldots,\vr_N$, yielding the expansion
\begin{equation}\label{eq:polesresidues}
\widetilde{\vu}(z)=\vs(z)+\sum_{k=1}^N\frac{\vr_k}{(z-\lambda)^k},
\end{equation}
for all $z$ in a punctured neighborhood of $\lambda$, with $\vs$ an analytic function. If $N=1$ (the ``simple'' case), the single corresponding residue of $\widetilde{\vu}$ can then be computed as
\begin{equation}\label{eq:residues}
\vr_1:=\lim_{z\to\lambda}(z-\lambda)\widetilde{\vu}(z)=\lim_{z\to\lambda}\frac{(z-\lambda)\vn(z)}{d(z)}=\frac{\vn(\lambda)}{d'(\lambda)}=-\sum_{j=1}^S\frac{q_j\vu(z_j)}{\lambda-z_j}\Bigg/\sum_{j=1}^S\frac{q_j}{(\lambda-z_j)^2}.
\end{equation}
If $N>1$, one first finds the highest-order residue $\vr_N$ using a formula similar to \cref{eq:residues}, and then recursively computes the other residues by deflation. We skip the details here.

Each pole-residue pair $(\lambda,\vr_k)$ satisfying \cref{eq:polesresidues} provides an approximation of an eigenpair of \cref{eq:ep}. See \cref{th:keldysh}.

\subsection{Adaptive sampling}\label{sec:21}
In the field of model reduction, it is customary to drive the addition of new sample points by an \emph{error indicator} $\rho:Z\to\xR$, with $Z\subset\xC$ being the region over which an approximation of the target $\vu$ is sought. More precisely, the following greedy sampling loop is performed:
\begin{center}
\begin{tabular}{ccccc}
initialize data set\hspace{-1em} & $\rightarrow$\hspace{-4.5em} & assemble surrogate $\widetilde{\vu}$ & $\rightarrow$ & assemble indicator $\rho$ \\[.5em]
 & \hspace{-1em} & $\uparrow$ & & $\downarrow$ \\[.5em]
 & \hspace{-1em} & add new sample $(z^\star,\vu(z^\star))$ to data set & $\leftarrow$ & find point $z^\star$ maximizing $\rho$ over $Z$\\
\end{tabular}
\end{center}

Here, the function $\vu$ that we approximate is the solution of the \emph{non-homogeneous} eigenproblem $\mT(z)\vu(z)=\vv$, cf. the definition of $\vu$ in \cref{sec:1}. For this kind of problems, it is common to employ the residual norm as indicator, i.e., $\rho(z)=\nnorm{\mT(z)\widetilde{\vu}(z)-\vv}_2$. However, its evaluation at a not-yet-seen point $z$ may, in general, be expensive, due to the need to assemble $\mT(z)$, to perform a matrix-vector product, and to compute a high-dimensional norm. Accordingly, the identification of the next sample point $z^\star$ may also be rather expensive.

To solve this issue, we rely on a different indicator. To motivate our choice of $\rho$, we first consider a simplified framework, namely, we assume that the target eigenproblem is \emph{linear}: $\mT(z)=\mT_0+z\mT_1$, with $\mT_0,\mT_1\in\xC^{n\times n}$. In such setting, it was shown in \cite{mri} that, as long as the surrogate $\widetilde{\vu}$ is an \emph{interpolatory} rational function $\widetilde{\vu}=\vn/d$ as in \cref{eq:approx}, the residual norm has the simplified form
\begin{equation}\label{eq:res}
\norm{\mT(z)\widetilde{\vu}(z)-\vv}_2=\frac{C(\mT_0,\mT_1,\vv,\vn)}{\aabs{d(z)}}=C(\mT_0,\mT_1,\vv,\vn)\Bigg/\Bigg|\sum_{j=1}^S\frac{q_j}{z-z_j}\Bigg|.
\end{equation}

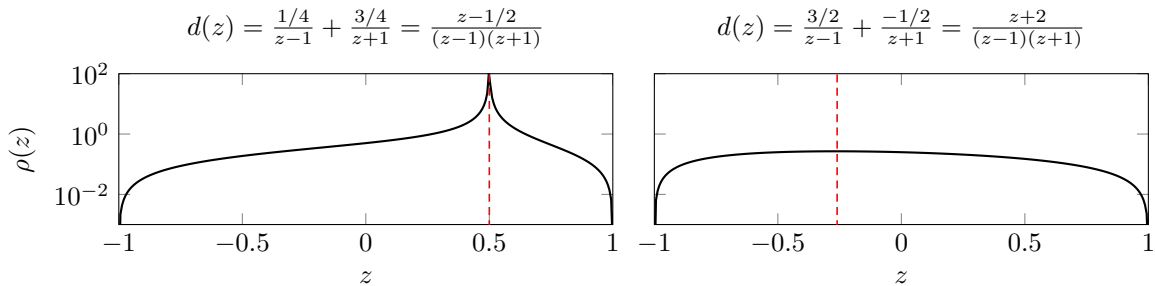
\begin{figure}[ht]
\pgfplotstableread[col sep=comma]{Data/test.csv}\tikzdata
\begin{tikzpicture}
\begin{semilogyaxis}[
	scale only axis,
	width = 6.5cm, height = 2cm,
	xlabel = {$z$},
	ylabel = {$\rho(z)$},
	xmin = -1, xmax = 1,
	ymin = 1e-3, ymax = 1e2,
	name = first plot,
	title = {$d(z)=\frac{1/4}{z-1}+\frac{3/4}{z+1}=\frac{z-1/2}{(z-1)(z+1)}$}
	]
\addplot[black, thick] plot table [x=z, y=r1]{\tikzdata};
\draw[red, densely dashed, semithick] (axis cs:.5, 1e-4) -- (axis cs:.5, 1e3);
\end{semilogyaxis}

\begin{semilogyaxis}[
	scale only axis,
	width = 6.5cm, height = 2cm,
	xlabel = {$z$},
	yticklabels = {,,},
	xmin = -1, xmax = 1,
	ymin = 1e-3, ymax = 1e2,
	at={(first plot.right of north east)},
	anchor=left of north west,
	name = second plot,
	title = {$d(z)=\frac{3/2}{z-1}+\frac{-1/2}{z+1}=\frac{\phantom{/2}z+2\phantom{/2}}{(z-1)(z+1)}$}
	]
\addplot[black, thick] plot table [x=z, y=r3]{\tikzdata};
\draw[red, densely dashed, semithick] (axis cs:-.26, 1e-4) -- (axis cs:-.26, 1e3);
\end{semilogyaxis}
%
\end{tikzpicture}
\caption{Simple rational residual indicator $\rho(z)=\aabs{d(z)}^{-1}$ in two test cases. The sample set consists of the points $\{-1,1\}$ in both plots. The greedy sampling loop, by maximizing the indicator, will place the next sample point at the location denoted by the vertical dashed line. Note that, in the left plot, the surrogate has a pole inside the interval $Z=[-1,1]$. The greedy indicator suggests taking a sample there.}
\label{fig:1}
\end{figure}

Since the constant $C$ in \cref{eq:res} does \emph{not} depend on $z$, one can find $z^\star=\argmax_{z\in Z}\nnorm{\mT(z)\widetilde{\vu}(z)-\vv}_2$ \emph{only by looking at the surrogate denominator $d$}. Taking a qualitative viewpoint, this way of finding the next sample point can be interpreted as follows: we add the next sample point by maximizing its distance from the already explored points (where information is already available) while minimizing its distance from the poles of the surrogates (where more information is needed to improve the accuracy of the eigenpair approximation). See \cref{fig:1} for an illustration.

If the problem is more general, namely, if $\mT(z)$ depends nonlinearly on $z$, then \cref{eq:res} is no longer true, since the residual norm may depend on $z$ in a more complicated way. However, we can still choose $\rho=\aabs{d(\cdot)}^{-1}$ as heuristic error indicator, and then use it to drive the greedy sampling loop.

\begin{remark}\normalfont
Our proposed greedy sampling driven by $\aabs{d(\cdot)}^{-1}$ has some affinity to Leja-Bagby (or similar) sampling schemes. See, e.g., \cite{lietaert}. However, the Leja-Bagby scheme chooses the next sample point by using information on the current sample set only. On the other hand, our approach, in the spirit of greedy model reduction, uses information on the current sample set \emph{but also on the current surrogate $\widetilde{\vu}$ itself}.
\end{remark}

\subsection{Theoretical considerations}\label{sec:22}
We devote this section to two topics. The first one is related to our claim that pole-residue pairs of $\vu(\cdot)=\mT(\cdot)^{-1}\vv$ provide information on eigenpairs of \cref{eq:ep}. We motivate this by the following particular version of Keldysh's theorem. See, e.g., \cite{beyn}.

\begin{theorem}\label{th:keldysh}
Let $\mT:\xC\to\xC^{n\times n}$ be analytic and assume that $\lambda\in\xC$ is an eigenvalue, i.e., there exists a solution of \cref{eq:ep} for such $\lambda$. For any $\vv\in\xC^n$, there exist an integer $N>0$, vectors $\{\vv_k\}_{k=1}^N\subset\xC^n\setminus\{\vzero\}$, coefficients $\{\alpha_k\}_{k=1}^N\subset\xC$, a neighborhood $Z$ of $\lambda$, and an analytic function $\vs:Z\to\xC^n$, such that
\begin{equation}\label{eq:keldysh}
\mT(z)^{-1}\vv=\vs(z)+\sum_{k=1}^N\frac{\alpha_k\vv_k}{(z-\lambda)^k}\quad\forall z\in Z\setminus\{\lambda\}.
\end{equation}
All vectors $\{\vv_k\}_{k=1}^N\subset\xC^n$ solve \cref{eq:ep} with eigenvalue $\lambda$.
\end{theorem}

From this result, we conclude that the eigenpair approximation described in \cref{sec:2} is sensible. Indeed, we are simply using surrogate quantities as approximations of the corresponding exact ones in \cref{eq:keldysh}.

To conclude this section, we would like to briefly discuss why we choose minimal rational interpolation for the construction of the approximation $\widetilde{\vu}$. Other interpolatory methods, like the (set-valued) Loewner framework \cite{mrib} or AAA \cite{aaa} use a barycentric form like \cref{eq:approx} but rely on least-squares formulations to find the coefficients $q_1,\ldots,q_S$. In this respect, minimal rational interpolation has the clear advantage of needing fewer sample points than its competitors, when employed to build a surrogate of certain degrees of $\vn$ and $d$: only $S$ samples are needed to approximate $S-1$ eigenpairs, whereas the other above-mentioned methods require approximately twice as many samples to achieve the same degree.

However, this efficiency comes at a price: from a theoretical viewpoint, the convergence results for minimal rational interpolation have restrictive assumptions. Indeed, the theory developed in \cite{mri} is limited to the case of simple poles ($N=1$ in \cref{th:keldysh}) with linearly independent residues. In particular, this implies that minimal rational interpolants might struggle with some kinds on nonlinearities. For instance, in a nonlinear eigenproblem of the form $\sin(\lambda)\vw_\lambda=\vzero$, minimal rational interpolation might not be able to identify more than a single eigenvalue, since all the eigenvectors of the problem are collinear. A more in-depth numerical investigation of these issues is in progress, and is outside the scope of the present paper.

\section{Numerical results}\label{sec:3}
In this section, we apply our proposed method to numerically solve a nonlinear eigeproblem stemming from an application in acoustics.

A plane pressure wave with wavenumber $k=10$ (all quantities are adimensional) impinges on the opening of a cavity, modeling a bottle-like Helmholtz resonator. The shape of the resonator is parametric, with the parameter $z\in Z=[1, 3]$ denoting the width of the neck of the resonator. The walls of the resonator are assumed to be sound-soft. We are interested in finding the values of $z$ (i.e., the geometric configurations) for which a resonance is observed. In layman's terms, we want to determine which resonator shapes ``make a sound'' in response to the given impinging wave.

\begin{figure}[t]
\hspace{1em}\includegraphics[scale=.3875]{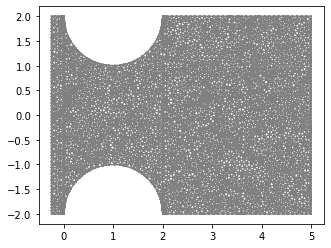}
\hspace{1em}\includegraphics[scale=.375]{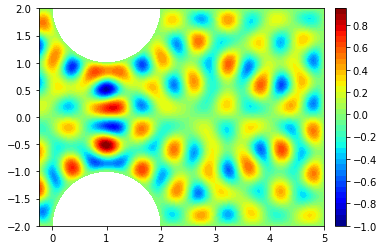}
\hspace{1em}\includegraphics[scale=.375]{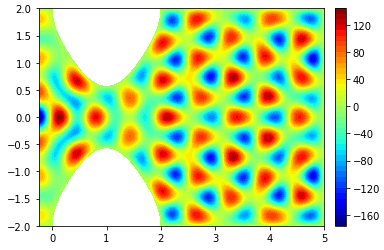}
\caption{Left: mesh on reference domain (coarsened to make the mesh visible). Center and right: solutions to the non-homogeneous problem for $z=2$ and $z=1.28$. The solution for $z=1.28$ has been remapped from $\Omega^\star$ to $\Omega(z=1.28)$ for illustration purposes. Note the scale.}
\label{fig:2}
\end{figure}

We base our mathematical model for this problem on the Helmholtz equation at wavenumber $k$, cast inside a 2D-section of the resonator $\Omega(z)\subset\xR^2$ (we assume that the problem is invariant with respect to the third space dimension). See \cref{fig:2} for a depiction of two sample geometric configurations. We specify our choice of parametrization for $\Omega$ in \cref{sec:a}. The boundary of the domain, namely, $\partial\Omega(z)$, is partitioned into inlet (where a Neumann boundary condition is imposed) and wall (where a Dirichlet boundary condition is imposed). The resulting problem reads
\begin{equation}\label{eq:pde}
\begin{cases}
(-\partial_{xx}-\partial_{yy}-k^2)w((x,y),z)=0&\text{for }(x,y)\in\Omega(z),\\
-\partial_x w((x,y),z)=0&\text{for }(x,y)\in\Gamma_{\text{inlet}}(z),\\
w((x,y),z)=0&\text{for }(x,y)\in\partial\Omega(z)\setminus\Gamma_{\text{inlet}}(z).
\end{cases}
\end{equation}
The solution field $w((x,y),z)$ is a function of the space coordinates $(x,y)\in\Omega(z)$, but also of the geometric parameter $z$. It represents the pressure intensity at a given point, for a given resonator geometry.

In order to handle the parametric geometry of the domain, we follow a mapping approach: we fix $\Omega^\star:=\Omega(2)$ as reference domain, introduce a nonlinear $z$-dependent bijective mapping $\phi_z:\Omega^\star\ni(\widehat{x},\widehat{y})\mapsto(x,y)\in\Omega(z)$, and cast \cref{eq:pde} in $\Omega^\star$. By setting $\widehat{w}((\widehat{x},\widehat{y}),z):=w(\phi_z(\widehat{x},\widehat{y}),z)$, we obtain
\begin{equation}\label{eq:pderef}
\begin{cases}
(-(a_1\partial_{\widehat{x}}+a_2\partial_{\widehat{x}})^2+(a_3\partial_{\widehat{x}}+a_4\partial_{\widehat{x}})^2-k^2)\widehat{w}((\widehat{x},\widehat{y}),z)=0&\text{for }(\widehat{x},\widehat{y})\in\Omega^\star,\\
-\partial_{\widehat{x}} \widehat{w}((\widehat{x},\widehat{y}),z)=0&\text{for }(\widehat{x},\widehat{y})\in\Gamma_{\text{inlet}}^\star:=\Gamma_{\text{inlet}}(2),\\
\widehat{w}((\widehat{x},\widehat{y}),z)=0&\text{for }(\widehat{x},\widehat{y})\in\partial\Omega^\star\setminus\Gamma_{\text{inlet}}^\star.
\end{cases}
\end{equation}
Note that the change of variables leads to the appearance of the coefficients $a_1,\ldots,a_4$, which are related to the Jacobian of $\phi_z$. In general, they may all be nonzero fields \emph{depending non-linearly on both $z$ and $(x,y)$}. The expressions of such coefficients are provided in \cref{sec:a}.

We introduce a triangulation of the domain $\Omega^\star$, and we approximate \cref{eq:pderef} using piecewise-quadratic finite elements. We employ FEniCS \cite{fenics} in our implementation. This turns the PDE into a linear system of the form $\mT\vw=\vzero$, where $\mT=\mT(z)\in\xC^{n\times n}$ is a sparse symmetric indefinite matrix with $z$-dependent entries and $\vw=\vw(z)\in\xC^n$ is a vector whose entries are nodal values of the finite-element approximation of $\widehat{w}$. For our discretization, $n=52344$. The task is now to find values of $z\in[1,3]$ for which the linear system has non-trivial solutions. Since $z$ enters $\mT$ through the nonlinear coefficients $a_1,\ldots,a_4$, this is a nonlinear eigenproblem like the one in \cref{eq:ep}.

In order to numerically approximate the eigenpairs, we introduce a vector $\vv\in\xC^n$, which we obtain by imposing non-homogeneous Neumann boundary conditions\footnote{Note that our choice of $\vv$ is motivated by the fact that the corresponding $\vu$ has physical meaning: it contains the nodal values of the finite-element solution of the Helmholtz equation corresponding to the imposed non-homogeneous boundary conditions. However, it is not necessary for the application of our technique: in fact, a random (e.g., Gaussian) vector would have served us perfectly well, if not better.}. Using the notation from \cref{sec:1}, we then define the vector-valued function $\vu(z)=\mT(z)^{-1}\vv$.

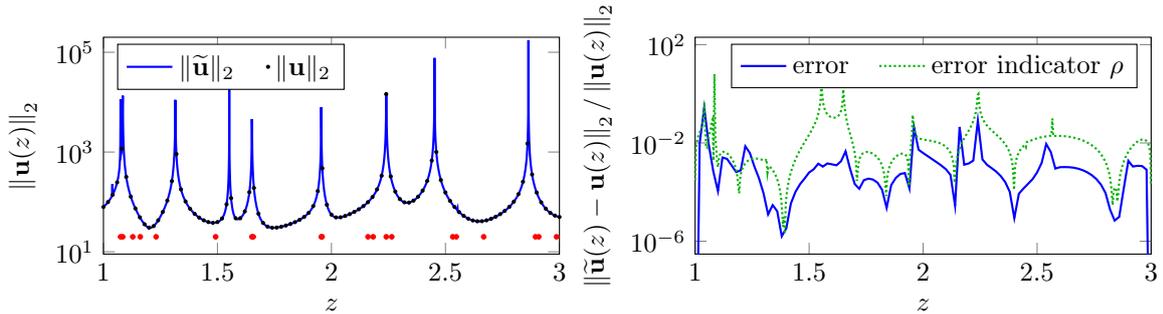
\begin{figure}[t]
\pgfplotstableread[col sep=comma]{Data/a_samples.csv}\tikzdatas
\pgfplotstableread[col sep=comma]{Data/a_error.csv}\tikzdata
\pgfplotstableread[col sep=comma]{Data/a_estimator.csv}\tikzdatae
\begin{tikzpicture}
\begin{semilogyaxis}[
	scale only axis,
	width = 6cm, height = 2.875cm,
	xlabel = {$z$},
	ylabel = {$\norm{\vu(z)}_2$},
	xmin = 1, xmax = 3,
	ymin = 9, ymax = 2e5,
	name = first plot,
	legend columns = 2,
	legend pos = north west,
	legend style={/tikz/every even column/.append style={column sep=1em}}
	]
\addplot[blue, thick] plot table [x=z, y=a20]{\tikzdata};
\addplot[black, only marks, mark=*, mark size=.2mm] plot table [x=zc, y=x]{\tikzdata};
\addplot[red, only marks, mark=*, mark size=.3mm] plot table [x=z40, y expr=2e1]{\tikzdatas};
\legend{$\nnorm{\widetilde{\vu}}_2$, $\nnorm{\vu}_2$};
\end{semilogyaxis}

\begin{semilogyaxis}[
	scale only axis,
	width = 6cm, height = 2.875cm,
	xlabel = {$z$},
	ylabel = {$\norm{\widetilde{\vu}(z)-\vu(z)}_2/\norm{\vu(z)}_2$},
	xmin = 1, xmax = 3,
	ymin = 3e-7, ymax = 2e2,
	at={(first plot.right of north east)},
	anchor=left of north west,
	legend columns = 2,
	legend pos = north east
	]
\addplot[blue, thick] plot table [x=zc, y=e20]{\tikzdata};
\addplot[green!70!black, densely dotted, thick] plot table [x=z20, y=e20]{\tikzdatae};
\legend{error \phantom{d}, error indicator $\rho$};
\end{semilogyaxis}
\end{tikzpicture}
\caption{Left plot: exact $\nnorm{\vu}_2$ (black dots) and approximated $\nnorm{\widetilde{\vu}}_2$ (blue line) using $20$ greedily selected samples. The $20$ sample points are denoted by the red dots at the bottom. Right plot: relative approximation error at $101$ validation points (blue line). The greedy indicator $\rho(z)=\aabs{d(z)}^{-1}$ (on a finer grid of 501 points) is also shown as a dotted line.}
\label{fig:3}
\end{figure}

We then run our proposed eigensolver. We select only the endpoints of $Z$, namely, $z_1=1$ and $z_2=3$, as initial sample set. We allocate a computational budget of $20$ samples, i.e., $20$ ``solves'' of the problem at different values of $z$, allowing the heuristic indicator from \cref{sec:21} to drive the selection of the sample points $z_3,\ldots,z_{20}$. We display the results of the greedy minimal rational approximation method in \cref{fig:3}. We observe a fairly good fit between $\vu$ and $\widetilde{\vu}$. Moreover, the heuristic estimator seems to capture rather well the locations where the approximation error is still large.

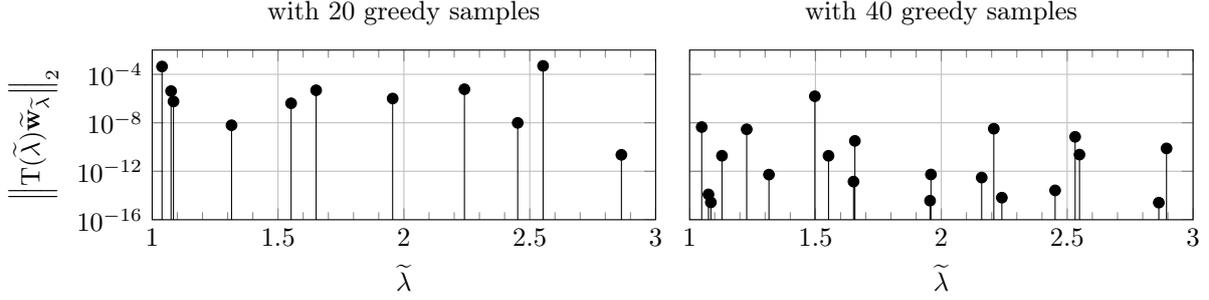
\begin{figure}[t]
\pgfplotstableread[col sep=comma]{Data/a_poles.csv}\tikzdata
\begin{tikzpicture}
\begin{semilogyaxis}[
	scale only axis, ycomb,
	minor x tick num=4,
	width = 6.625cm, height = 2.25cm,
	xlabel = {$\widetilde{\lambda}$},
	ylabel = {$\norm{\mT(\widetilde{\lambda})\widetilde{\vw}_{\widetilde{\lambda}}}_2$},
	title = {with $20$ greedy samples},
	xmin = 1, xmax = 3,
	ymin = 1e-16, ymax = 1e-2,
	ytick = {1e-16,1e-12,1e-8,1e-4},
	grid=major,
	name = first plot, log origin=infty
	]
\addplot[black, mark=*] plot table [x=p20, y=r20]{\tikzdata};
\end{semilogyaxis}

\begin{semilogyaxis}[
	scale only axis, ycomb,
	minor x tick num=4,
	width = 6.625cm, height = 2.25cm,
	xlabel = {$\widetilde{\lambda}$},
	ylabel = {},
	title = {with $40$ greedy samples},
	yticklabels = {,,},
	xmin = 1, xmax = 3,
	ymin = 1e-16, ymax = 1e-2,
	ytick = {1e-16,1e-12,1e-8,1e-4},
	grid=major,
	at={(first plot.right of north east)},
	anchor=left of north west,
    log origin=infty
	]
\addplot[black, mark=*] plot table [x=p40, y=r40]{\tikzdata};
\end{semilogyaxis}
\end{tikzpicture}
\caption{Magnitude of residual of homogeneous eigenproblem. Each pin is placed at some approximated eigenvalue $\widetilde{\lambda}$, with $\widetilde{\vw}_{\widetilde{\lambda}}$ being the corresponding (normalized) eigenvector. In the left and right plots we show the results obtained with $20$ and $40$ greedy samples, respectively.}
\label{fig:4}
\end{figure}

Concerning the approximation of the eigenpairs, our method identifies $11$ eigenvalues within the region $Z$, as shown in \cref{fig:4} (left). In addition, $8$ more eigenvalues are identified outside $Z$, but we do not show them here. A comparison of \cref{fig:3} (left) and \cref{fig:4} (left) shows that the sample points are (loosely) placed around the (\emph{a priori} unknown) eigenvalues.

As a way to assess the quality of the approximation, we evaluate the norm of the residual of the \emph{homogeneous} eigenproblem, namely, $\nnorm{\mT(\widetilde{\lambda})\widetilde{\vw}_{\widetilde{\lambda}}}_2$, where $(\widetilde{\lambda},\widetilde{\vw}_{\widetilde{\lambda}})$ is an approximated eigenpair (with normalized eigenvector), obtained from a pole-residue pair of the rational surrogate $\widetilde{\vu}$. In \cref{fig:4} (left) we display the magnitude of such residuals. We can observe that the residual norms range between $10^{-11}$ and $10^{-3}$, in a way that is somewhat consistent with the error curve in \cref{fig:3}. This is perhaps slightly surprising, since the latter curve depicts the approximation error for the \emph{non-homogeneous} problem.

\begin{figure}[t]
\pgfplotstableread[col sep=comma]{Data/a_samples.csv}\tikzdatas
\pgfplotstableread[col sep=comma]{Data/a_error.csv}\tikzdata
\pgfplotstableread[col sep=comma]{Data/a_estimator.csv}\tikzdatae
\begin{tikzpicture}
\begin{semilogyaxis}[
	scale only axis,
	width = 6cm, height = 2.875cm,
	xlabel = {$z$},
	ylabel = {$\norm{\vu(z)}_2$},
	xmin = 1, xmax = 3,
	ymin = 9, ymax = 2e5,
	name = first plot,
	legend columns = 2,
	legend pos = north west,
	legend style={/tikz/every even column/.append style={column sep=1em}}
	]
\addplot[blue, thick] plot table [x=z, y=a40]{\tikzdata};
\addplot[black, only marks, mark=*, mark size=.2mm] plot table [x=zc, y=x]{\tikzdata};
\addplot[red, only marks, mark=*, mark size=.3mm] plot table [x=z20, y expr=2e1]{\tikzdatas};
\addplot[red, only marks, mark=*, mark size=.3mm] plot table [x=z40, y expr=2e1]{\tikzdatas};
\legend{$\nnorm{\widetilde{\vu}}_2$, $\nnorm{\vu}_2$};
\end{semilogyaxis}

\begin{semilogyaxis}[
	scale only axis,
	width = 6cm, height = 2.875cm,
	xlabel = {$z$},
	ylabel = {$\norm{\widetilde{\vu}(z)-\vu(z)}_2/\norm{\vu(z)}_2$},
	xmin = 1, xmax = 3,
	ymin = 3e-13, ymax = 2e-4,
	at={(first plot.right of north east)},
	anchor=left of north west,
	legend columns = 2,
	legend pos = north east
	]
\addplot[blue, thick] plot table [x=zc, y=e40]{\tikzdata};
\addplot[green!70!black, densely dotted, thick] plot table [x=z40, y expr=\thisrow{e40}*5]{\tikzdatae};
\legend{error \phantom{d}, error indicator $\rho$};
\end{semilogyaxis}
\end{tikzpicture}
\caption{Left plot: exact $\nnorm{\vu}_2$ (black dots) and approximated $\nnorm{\widetilde{\vu}}_2$ (blue line) using $40$ greedily selected samples. The $40$ sample points are denoted by the red dots at the bottom. Right plot: relative approximation error at $101$ validation points (blue line). The greedy indicator $\rho(z)=\aabs{d(z)}^{-1}$ (on a finer grid of 501 points) is also shown as a dotted line.}
\label{fig:3b}
\end{figure}
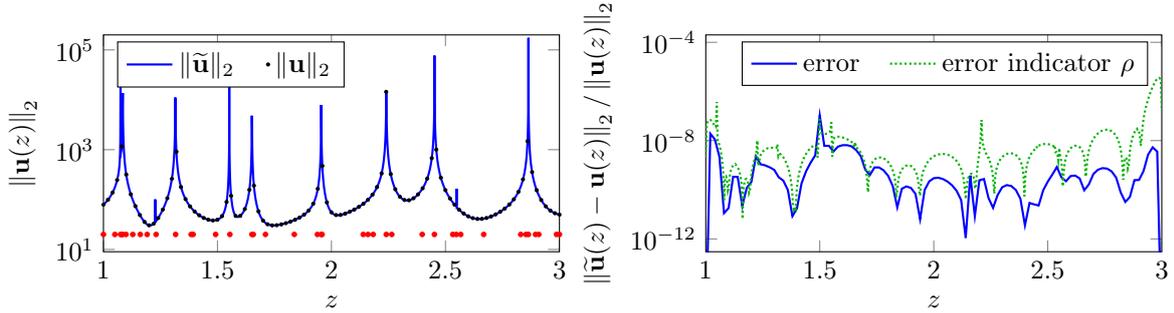

As an additional numerical test, we double the computational budget and repeat our experiment. We show the corresponding approximation results in \cref{fig:3b}. We note that the approximation error (for the non-homogeneous problem) decreases by about 5 orders of magnitude. In \cref{fig:4} (right) we show the results pertaining to the eigenpairs and to the residual of the homogeneous eigenproblem. We observe two effects: an overall decrease of the norm of the residual (again, by around 5 orders of magnitude) and the appearance of (seemingly) spurious eigenpairs. The latter is a rather serious issue, which could be solved (or, at least, weakened) by introducing a ``filtering'' step in post-processing, with the aim of eliminating multiple approximations of the same eigenvalue.

\section{Conclusions}\label{sec:4}
We have presented a technique for the numerical approximation of solutions of generic nonlinear eigenproblems. Being based on the rational approximation of the solution of a non-homogeneous version of the eigenproblem, our method is \emph{non-intrusive}. In particular, it can be applied regardless of \emph{how} $\mT(z)$ depends on $z$. Despite this, our method is equipped with a heuristic indicator that allows for an adaptive construction of the sample set. In our numerical test, we have showcased the effectiveness of our method and of our indicator. However, we have also observed the appearance of spurious effects if the number of samples is disproportionate to the (unknown) number of eigenvalues to be approximated. Possible solutions to this issue might borrow some ideas from state-of-the-art approaches for automatic degree reduction of rational approximants. See, e.g., \cite{rkfit}. We leave this as a direction for future research.

\appendix

\section{Appendix}\label{sec:a}
We fix the reference domain $\Omega^\star=[-\frac14,5]\times[-2,2]\setminus B(1,2)\setminus B(1,-2)$, with $B(\widehat{x},\widehat{y})\subset\xR^2$ the unit disk centered at $(\widehat{x},\widehat{y})$. The corresponding inlet boundary is $\Gamma_{\text{inlet}}^\star=\{-\frac14\}\times[-2,2]$. For any $z\in[1,3]$, we define the $C^1(\Omega^\star)$-regular mapping
\begin{equation*}
\phi_z((\widehat{x},\widehat{y}))=\begin{cases}
\left(\widehat{x},\left(\frac{2+z}4+\frac{2-z}4\cos(\pi\widehat{x})\right)\widehat{y}\right)&\text{if }0<\widehat{x}<2,\\
(\widehat{x},\widehat{y})&\text{otherwise},
\end{cases}
\end{equation*}
and the deformed domain $\Omega(z)=\phi_z(\Omega^\star)$. The Jacobian of $\phi_z$ is the identity for $\widehat{x}\leq 0$ or $\widehat{x}\geq 2$, whereas
\begin{equation*}
J\phi_z((\widehat{x},\widehat{y}))=\begin{bmatrix}
1 & 0\\
\pi\frac{z-2}4\sin(\pi\widehat{x})\widehat{y} & \frac{2+z}4+\frac{2-z}4\cos(\pi\widehat{x})
\end{bmatrix}\quad\text{if }0<\widehat{x}<2.
\end{equation*}
Since $[\partial_x,\partial_y]=[\partial_{\widehat{x}},\partial_{\widehat{y}}]\left(J\phi_z((\widehat{x},\widehat{y}))\right)^{-1}$ by the chain rule, we set
\begin{equation*}
\begin{bmatrix}
a_1((\widehat{x},\widehat{y}),z) & a_2((\widehat{x},\widehat{y}),z)\\
a_3((\widehat{x},\widehat{y}),z) & a_4((\widehat{x},\widehat{y}),z)
\end{bmatrix}:=\left(J\phi_z((\widehat{x},\widehat{y}))\right)^{-\top}.
\end{equation*}

This gives
\begin{equation*}
a_1((\widehat{x},\widehat{y}),z)=1,\ a_3((\widehat{x},\widehat{y}),z)=0,\ a_4((\widehat{x},\widehat{y}),z)=\begin{cases}
\left(\frac{2+z}4+\frac{2-z}4\cos(\pi\widehat{x})\right)^{-1}&\text{if }0<\widehat{x}<2,\\
1&\text{otherwise},
\end{cases}
\end{equation*}
and
\begin{equation*}
a_2((\widehat{x},\widehat{y}),z)=\begin{cases}
\left(\pi\frac{2-z}4\sin(\pi\widehat{x})\widehat{y}\right)\big/\left(\frac{2+z}4+\frac{2-z}4\cos(\pi\widehat{x})\right) &\text{if }0<\widehat{x}<2,\\
0&\text{otherwise}.
\end{cases}
\end{equation*}

\vspace{\baselineskip}

\bibliographystyle{plain}
\bibliography{paper.bib}

\end{document}